\documentclass[11pt]{amsart}
\usepackage[utf8]{inputenc}
\usepackage{amsmath,amssymb,amsfonts,anysize,mathrsfs,graphicx,setspace,parskip,multicol,upgreek,yfonts,tikz,physics,graphicx,xcolor,csquotes,enumerate,indentfirst,ifpdf}
\usepackage[bottom]{footmisc}
\pdfoutput=1
\newtheorem{theorem}{Theorem}
\newtheorem{lemma}[theorem]{Lemma}

\newtheorem{proposition}[theorem]{Proposition}
\newtheorem{question}{Question}
\newtheorem{remark}{Remark}

\newtheorem{definition}{Definition}
\newtheorem{fact}{Fact}
\newtheorem{claim}{Claim}

\usepackage{hyperref}
\hypersetup{colorlinks=true, citecolor=blue}

\title{$Q$-sets, $\Delta$-sets and $L$-spaces}

\author{P. Memarpanahi}
\address{Department of Mathematics and Statistics, York University, Toronto, Ontario M3J 1P3, Canada}

\author{P. Szeptycki}
\address{Department of Mathematics and Statistics, York University, Toronto, Ontario M3J 1P3, Canada}
\email{szeptyck@yorku.ca}
\thanks{2020 \emph{Mathematics Subject Classification}. Primary 54D20, 54H05, 
 03E35; Secondary 03E75, 54E52, 54G20. }
\thanks{The second author acknowledges support from NSERC grant 503970.}

\begin{document}

\maketitle

\begin{abstract} The question whether there is a Lindelof Q-set space or Lindelof $\Delta$-set space is considered. We show that J. Moore's ZFC $L$-space is not a Q-set space in ZFC and, assuming all Aronszajn trees are special, it is not a $\Delta$-set space. 
\end{abstract}

\section{Introduction}

An uncountable subset $A$ of $\mathbb R$ is called a {$Q$-set} if every subset of $A$ is a relative $G_\delta$. In general, an uncountable regular topological space $X$ is called a {$Q$-set space} if $X$ is not $\sigma$-discrete and every subset of $X$ is a $G_\delta$ set. 

The existence of a $Q$-set of size $\kappa$ implies that $2^{\kappa}=2^{\omega}$, so it is consistent, e.g., assuming the Continuum Hypothesis (CH), that there are no $Q$-sets. 
 Under Martin Axiom (MA) and the negation of CH, every uncountable subset $A$ of $\mathbb R$ of size smaller than $\mathfrak c$ is a $Q$-set (see \cite{M}). Much later, Zolt\'{a}n Balogh employed his powerful elementary submodal technique, also used to construct a variety of small ZFC Dowker spaces \cite{Balogh-Dowker}, to construct a $\operatorname{ZFC}$ example of a $Q$-set space \cite{Balogh-Q-set}. Later, he improved this construction to obtain a paracompact $Q$-set space \cite{Balogh-paracompact}, and asked whether there is a Lindel\"of $Q$-set space in ZFC \cite{B-note}. This was also asked by Leidermann and Tkacenko in \cite{AT}:
 
\begin{question} (Z. Balogh) Does there exist a Lindel\"of $Q$-set space in $\operatorname{ZFC}$? 
    
\end{question} 

Note that a $Q$-set of reals would be a Lindel\"of $Q$-set space. Moreover, any Lindel\"of Q-set space would necessarily be hereditarily Lindel\"of. While it is consistent with CH that there is a separable $Q$-set space, V=L implies there are none (see proposition 3.1 in \cite{Balogh-Junilla}) and, moreover, hereditarily separable Q-set spaces must have cardinality $<{\mathfrak c}$ (\cite{B-note}).  Thus, a Lindelof Q-set space constructed assuming only ZFC could not be separable and, therefore, would be an $L$-space. 

  While there are quite a number of consistent examples of L-spaces, the celebrated example of J. Moore \cite{L-space} is essentially the only known ZFC example. In the first section, we will show that Moore's L-space is not a $Q$-set space, and so, providing an answer to Question 1 would require the construction of an essentially different example of an L-space. 
  
Some might first search for a consistent example of an L-space that is also a Q-set space; however, Balogh proved the following in his unpublished notes. 
\begin{proposition}\cite{B-note}
    If there is an L-space in $\operatorname{ZFC}$, then under $\operatorname{MA}(\omega_1),$ there exists a Lindel\"{o}f $Q$-set space which is not hereditarily separable.  
    \begin{proof}
Let $(X, \tau_0)$ be an L-space that is not hereditarily separable. Without loss of generality, assume that $X$ is not separable. Furthermore, suppose $\abs{X} = \aleph_1$; hence, it is zero-dimensional, being a completely regular (regular + Lindel\"{o}f) space. 

Now, take a subset A of the reals of size $\aleph_1$, and identify it with $X$. This gives us a zero-dimensional, separable, metrizable topology $(X, \tau_E)$ inherited from the Euclidean topology on $\mathbb{R}$.

Let $(X, \tau_\cup)$ be the topology generated by the union of $\tau_0$ and $\tau_E$. Clearly, $(X, \tau_\cup)$ is not separable, but is zero dimensional since every $U\in\tau_\cup$ is some union of finite intersection of elements of zero dimensional topologies $\tau_0$ and $\tau_E.$ 

Furthermore, $X$ being a $Q$-set with respect to $\tau_E$ would imply that every of its subset is an $F_\sigma$ set in $(X,\tau_E).$ Now since closed subsets of the space $(X,\tau_E)$ remain closed in $\tau_\cup,$ it follows that $X$ is a $Q$-set with respect to $\tau_\cup.$
It remains to show that $(X,\tau_\cup)$ is hereditarily Lindel\"{o}f. To this end, fix a countable basis ${\mathcal B}$ for $\tau_E$ and let $\mathcal U=\{U_\alpha:\alpha\in J\}$ be a collection of basic open sets $U_\alpha$ in $\tau_\cup$. For each $\alpha\in J$ there are $E^\alpha\in {\mathcal B}$ and $O^\alpha\in \tau_0$ such that
$$U_\alpha=E^\alpha\cap O^\alpha$$
For each $E\in {\mathcal B}$, let ${\mathcal V}_E=\{O^\alpha \in \tau_0:E=E^\alpha\}$. Since $\tau_0$ is a hereditarily Lindelof topology, there is a countable family of sets ${\mathcal O}_E\subseteq {\mathcal V}_E$ such that $\bigcup {\mathcal O}_E=\bigcup{\mathcal V}_E$. Since ${\mathcal B}$ is countable, $${\mathcal W}=\{E\cap O:E\in{\mathcal B}\text{ and }O\in {\mathcal O}_E\}$$
is a countable subset of ${\mathcal U}$ with $\bigcup {\mathcal W}=\bigcup{\mathcal U}$.
\end{proof}
\end{proposition}

Closely related to the concept of a $Q$-set is that of a $\Delta$-set \cite{Delta-paracompact} and, more generally, a $\Delta$-space \cite{AT}.

\begin{definition} A topological space $X$ is a $\Delta$-space if for every decreasing sequence of subsets $\{D_n:~n\in\omega\}$ of $X$, with empty intersections, there exists a decreasing sequence of open sets $\{U_n:~U_n\supseteq D_n\}$ of $X$ such that  
\[\displaystyle\bigcap_{n\in\omega}U_n=\emptyset.\]
A $\Delta$-set is a subspace of the real line that is a $\Delta$-space.
\end{definition}

It is not hard to see that every $Q$-set space is a $\Delta$-space \cite{Q-set-delta-set}. And while there are several examples of $\Delta$-spaces in ZFC, they all fail to be $Q$-set spaces since they are either $\sigma$-discrete or they have points or closed subsets that are not $G_\delta$. So, we will call a nontrivial $\Delta$-space (i.e., not $\sigma$-discrete and with closed subsets $G_\delta$) a {\em $\Delta$-set space}. The notion of a $\Delta$-set of reals arose in the study of the normal Moore space conjecture, where $Q$-sets were used to construct important counterexamples to the conjecture. For example, the tangent disc topology over a set $X$ is normal if and only if $X$ is a $Q$-set, while the space is countably paracompact if and only if $X$ is a $\Delta$-set (see \cite{Delta-paracompact}). More recently, the class of $\Delta$-spaces has gained interest in the study of $C_p$-spaces, as the class of all $\Delta$-spaces consists precisely of those $X$ for which the locally convex space $C_p(X)$ is distinguished \cite{KL}.

Let us first point out that many of the consistent examples of $L$-spaces are not even $\Delta$-set spaces. Recall that a space is said to be resolvable if it can be partitioned into two dense subsets, and $\omega$-resolvable if the space admits a countable pairwise disjoint family of dense subsets. 
Filatova showed that a Lindel\"{o}f space of uncountable dispersion character is resolvable \cite{resolvable} and this was improved to $\omega$-resolvable in \cite{omega-resolvable}. It is easy to see -- and was proved in \cite{AP} -- that any Baire Lindel\"{o}f space without an isolated point is not a $\Delta$-space. One can easily check that many consistent examples of $L$-spaces are clearly Baire spaces (e.g., Luzin spaces and other examples presented in \cite{classical-L-spaces}) and, therefore, not $\Delta$ or $Q$.

One exception is the HFC $L$ subspaces of $2^{\omega_1}$  due to Juh\'{a}sz (see \cite{classical-L-spaces}). While it is easy to construct a Baire HFC assuming CH, we don't know whether every HFC is Baire. 

\begin{question} Is there an HFC that is not a Baire space. Can an HFC be a Q-set space or $\Delta$-set space?
\end{question}

The following is a natural weakening of Questions 1 and 2 and was raised in \cite{AP}
\begin{question}
    Does there exist a Lindel{\"o}f $\Delta$-set space in ZFC? Or even a consistent example of an L-space that is also a $\Delta$-space?
\end{question}
Or perhaps they are equivalent questions: 

 \begin{question} If $X$ is Lindel\"of, is $X$ being a $Q$-set space equivalent to $X$ being a $\Delta$-set space?
 \end{question}

Of course, this relates to the long-standing question posed by Reed \cite{Delta-paracompact} whether there exists a $\Delta$-set that is not a $Q$-set. Although a consistent counterexample has appeared in the literature, many details of the construction seem unclear. Therefore, we believe that this question is still open, even for sets of reals. 

 Our other main result is that it is consistent that Moore's $L$-space is not a $\Delta$-set space (e.g., under the assumption that all Aronszajn trees are special), but we do not know whether, for example in the model obtained by adding one Cohen real, it could be a $\Delta$-set space. Indeed, we show that in this model, the Aronszajn tree naturally associated with Moore's L-space is not special. 

 \section{Moore's L-space}
Before we turn to the proof that Moore's L-space is not a Q-set space, we recall some fundamental notions used in the construction of the space. The reader interested in only this specific result may skip over these details, as they are primarily needed for results in later sections. 

First, we will recall the notion of a minimal walk between ordinals $\alpha, \beta\in \omega_1$ (see \cite{L-space} or \cite{minimal-walk}).
A walk from the ordinal $\alpha$ to the ordinal $\beta<\alpha$ is simply a finite sequence of ordinals $\langle \beta_i:0\leq i\leq n \rangle$ with $\beta=\beta_0>\cdots>\beta_n=\alpha$. We wish to consider those walks defined from a $C$-sequence:
\begin{definition} A sequence $\langle C_\alpha:\alpha\in\omega_1\rangle$ is called a {\em $C$-sequence} if given any $\alpha$, $C_\alpha$ is a cofinal subset of $\alpha$ such that: 
\begin{enumerate}
    \item\text{If $\alpha$ is a limit ordinal, then $C_\alpha\cap\xi$ is finite for all $\xi<\alpha$ }
    \item\text{If $\alpha=\beta+1$, then $C_\alpha=\{\beta\}$}
\end{enumerate}
\end{definition}

\begin{definition}[Upper Trace Function]\cite{minimal-walk}  It is the function $U:[\omega_1]^2\to[\omega_1]^{<\omega}$ such that for all $\alpha\leq\beta$:
\begin{align*}
&U(\alpha,\alpha)=\emptyset\\
&U(\alpha,\beta)=U(\alpha,\min{(C_\beta\setminus\alpha)})\cup\{\beta\}
\end{align*}
    
\end{definition}

\begin{definition}[Lower Trace Function] \cite{L-space} It is the function $L:[\omega_1]^2\to[\omega_1]^{<\omega}$ such that for all $\alpha\leq\beta:$
\begin{align*}
    &L(\alpha,\alpha)=\emptyset\\
    &L(\alpha,\beta)=L(\alpha,\min{(C_\beta\setminus\alpha)})\cup\big(\{\max(C_\beta\cap\alpha)\}\setminus{\max(C_\beta\cap\alpha)}\big).
\end{align*}
\end{definition}
\begin{remark} Note that the $k^{th}$ term $\xi_k$ of the lower trace $L(\alpha,\beta)$ could have also been defined as $\xi_k=\max(\cup_{i=0}^k(C_{\beta_i}\cap\alpha))$, where $\beta_i\in U(\alpha,\beta)$. Furthermore, this is an increasing function whereas the upper trace function is a strictly decreasing function.  
    
\end{remark}

The following two facts can be easily established by employing the definition of the lower trace function (see \cite{L-space} for further details).

\begin{fact} Given $\alpha\leq\gamma\leq\beta$, if $L(\gamma,\beta)<L(\alpha,\gamma)$ then $L(\alpha,\beta)=L(\alpha,\gamma)\cup L(\gamma,\beta)$.
\end{fact}
\begin{fact}\label{fact2}
    Let $0<\beta$ be a limit ordinal and $\alpha<\beta$ such that $L(\alpha,\beta)\neq\emptyset$. Then we have that $\displaystyle\lim_{\alpha\to\beta}\min L(\alpha,\beta)=\beta.$
\end{fact}
\begin{remark} An important consequence of Fact 2, which will be used often, is that if $\beta$ is a limit ordinal and  $\gamma>\beta$, we can then find an ordinal $\alpha$ such that $L(\beta,\gamma)<L(\xi,\beta)$ $\forall~\alpha\leq\xi<\beta.$
    
\end{remark}
The functions $e(\alpha,\beta)=e_\beta(\alpha)$, defined in the following definition, is referred to as \textquote{maximum weight functions} (see \cite{coherent}).

\begin{definition}\label{max} Given a $C$-sequence $\langle C_\alpha:\alpha\in\omega_1\rangle$ and $\beta\in\omega_1$ let $e_\beta:\alpha\rightarrow \omega$ be defined by $e_\beta(\alpha)=\max\{\abs{C_\xi\cap\alpha}:\xi\in U(\alpha,\beta)\}$.
\end{definition}

\begin{fact}\cite{coherent} These functions are finite-to-one and form a coherent family of functions, where coherent means that the set $\{\alpha\in\gamma:e_\gamma(\alpha)\neq e_\beta(\alpha)\}$ is finite for any pair of ordinals $\gamma<\beta$.
\end{fact}

\vskip2mm

Lastly, we are going to introduce two important functions that are directly used in constructing Moore's $L$-space.
\begin{definition}[Oscillation Function]\cite{L-space} Let $\langle e_\beta:\beta\in\omega_1\rangle$ be a sequence of functions as defined above. Define the oscillation function $\operatorname{osc}:[\omega_1]^2\rightarrow\omega$, to be the map so that for every $\alpha<\beta$:

$$\operatorname{osc}(\alpha,\beta)=\abs{\{\zeta\in{L(\alpha,\beta)\setminus\{\min L(\alpha,\beta)\}}:e_\alpha(\zeta^{L-})\leq{e_\beta(\zeta^{L-})}\wedge{e_\alpha}(\zeta)>e_\beta(\zeta)\}}\footnote{\text{We refer to $\operatorname{osc}(\alpha,\beta)$ explicitly as the number of oscillations between the pair functions $(e_\alpha,e_\beta)$ on the lower trace $L(\alpha,\beta).$}}$$

\end{definition}
Recall that for $\alpha<\beta$ and $\zeta\in L(\alpha,\beta)$ we used $\zeta^{L-}$ to denote the immediate predecessor of $\zeta$ in the corresponding lower trace.

Fix an uncountable rationally independent $\{z_\alpha:\alpha<\omega_1\}$ sequence in $\mathcal{S}$ the unit circle of the complex plane. Now we introduce the following function that defines Moore's $L$-space: 
\begin{definition}[o-Function]\cite{L-space} The function $o:[\omega_1]^2\to\mathcal S=\{z\in\mathbb C:z\overline z=1\}$ is defined as follow: 
$$o(\alpha,\beta)=z_\alpha^{\operatorname{osc}(\alpha,\beta)+1}~\text{where $\alpha\leq\beta$ }$$

\end{definition}
\vskip 6pt

Furthermore, define $w_\beta(\alpha)=\begin{cases}
o(\alpha,\beta),~\text{if $\alpha<\beta$}
\\1,\quad\text{otherwise}
    
\end{cases}$ 
\vskip 6pt

Now Moore's L-space is the space $$\mathscr L=\{w_\beta : \beta\in\omega_1\}\subseteq {\mathcal S}^{\omega_1}$$

We will invoke the following important proposition and theorem to prove our first result.
\begin{proposition}\cite{L-space}
    $T(o)=\{w_\beta\restriction\gamma:\gamma\leq\beta<\omega_1\}$ is an Aronszajn tree.
\end{proposition}

\begin{theorem}\label{thm3}\cite{L-space} Let $E$ be an uncountable subset of $\mathscr L$. Then the closure of $E$ in $\mathcal S^{\omega_1}$ contains a
set of the form $\{f\in\mathcal S^{\omega_1} : f\upharpoonright\alpha = \tau\}$ for some $\alpha<\omega_1$ and $\tau\in\mathcal S^\alpha$. Moreover, $\tau$ can be chosen so that $\{f\in E : f\upharpoonright\alpha = \tau\}$ is uncountable.
\end{theorem}

As stated, Theorem \ref{thm3} is formally stronger than the analogous result in \cite{L-space}, but Moore's proof does give the stronger statement. For completeness sake we will give a sketch of the proof. 

Let us first recall an important claim used in the proof. 
 \begin{claim} Let $T(E)$ be the set of all $\tau$ in $T(o)$ such
that for uncountably many $w_\beta$ in $E$, $\tau$ is a restriction of $w_\beta$, and for every element $\tau$ of $T(E)$, let $E_\tau=\{w_\beta\in E:~\tau=w_\beta\upharpoonright\alpha\}$, where $\alpha$ is the height of $\tau$. Then for each $\tau$ in $T(E)$ there is an $\alpha_\tau<\omega_1$ such that the projection map
$\pi:\mathcal S^{\omega_1}\to\mathcal S^{\omega_\setminus\alpha_\tau}$ sends the closure of $E_\tau$ onto $\mathcal S^{\omega_1\setminus\alpha_\tau}$ (see claim 7.16 \cite{L-space} for a proof).
\end{claim}

Now let us prove Theorem \ref{thm3}.

 \begin{proof}  For every $\tau$ in $T(E),$ there is an $\alpha_\tau$ by Claim $1$. Let $f$ be the function that sends each $\tau$ to the minimum $\alpha_\tau$. Consider the set of all $\tau\in T(E)$ of height less than $\omega_1$. Since the levels of $T(o)$ are countable, the closure of this set under this function $f$ will be countable. Hence, we can find a $\delta<\omega_1$ such that if  $\tau$ is in $T(E)$ of height less than $\delta,$ then the $\alpha_\tau$ obtained by Claim 1 will also be less than $\delta$. 
 
We fix $\tau_0\in T(E)$ of height $\delta$ with uncountably many extensions in $E$, and we now show that any $f\in\mathcal S^{\omega_1}$ that extends $\tau_0$ is in the full closure $cl(E)$ of $E$. 
 
 To see this, fix any $f$ extending $\tau_0$ and fix $\xi<\delta$. By choice of $\delta,$ we have  $\alpha_{\tau_0\upharpoonright\xi}<\delta$ and, by the properties of the $\alpha_\tau$'s, we have that the projection of $cl(E_{\tau_0\upharpoonright\xi})$ by the map $\pi_{\alpha_{\tau_0\upharpoonright\xi}}$ is onto $S^{\omega_\setminus\alpha_{\tau_0\upharpoonright\xi}}$. Therefore, for any such $\xi,$ there is a $g_\xi$ that agrees with $f$ both below $\xi$ and above $\alpha_{\tau_0\upharpoonright\xi}$, which means that $g_\xi$ is in the closure of $E$. So, we recursively choose $\xi_n$ increasing so that $\xi_{n+1}$ is above 
 $\alpha_{\tau_0\upharpoonright\xi_n}$. Then the sequence of $g_{\xi_n}$ converges to $f$ and since each $g_{\xi_n}$ is in the closure of $E$, it follows that $f$ is in the closure of $E$, completing the proof. 
  \end{proof}

\begin{theorem} $\mathscr L$ is not a $Q$-set space.
\end{theorem}
\begin{proof}  Suppose $\mathscr L$ is a $Q$-set space. Let $\mathscr{F}=\{A_\gamma:\gamma<\omega_2\}$ be an almost disjoint family of uncountable subsets of $\mathscr{L}$, i.e., so that for every $\gamma\neq\alpha$ in $ \omega_2$, we have $\abs{A_\gamma\cap A_\alpha}\leq\aleph_0.$ Since every subset in $\mathscr L$ is an $F_\sigma$, we can, without loss of generality, assume that each $A_\gamma$ is closed. Theorem \ref{thm3} guarantees that for every $\gamma<\omega_2,$ there is a pair $(\alpha_\gamma,\tau_\gamma)$ such that $\alpha_\gamma\in \omega_1$ and $\tau_\gamma\in \mathcal S^{\alpha_\gamma}$ such that the closure of $A_\gamma$ in $\mathcal S^{\omega_1}$ contains the set $\{f\in\mathcal{S}^{\omega_1}:f\restriction\alpha_\gamma=\tau_\gamma\}$. Moreover, for each $\gamma,$ we can choose an $\tau_\gamma\in T(o)$ in such a way that uncountably many elements of $A_\gamma$ extend $\tau_\gamma$. Since $\abs{T(o)}= \aleph_1,$ we can find a subset $J$ of $\omega_2$ of cardinality $\aleph_2$ and a $\tau\in T(o)$ of some height $\delta<\omega_1$ such that for every $\gamma\in J,$ we have that $\tau=\tau_\gamma$. Consequently, the closure of $A_\gamma$ in $\mathcal S^{\omega_1}$ contains the set $D=\{f\in\mathcal{S}^{\omega_1}:f\restriction\delta=\tau\}$, which again as mentioned, contains uncountably many elements of $A_\gamma$ for every $\gamma\in J$. Furthermore, due to the fact that the levels of $T(o)=\{w_\beta\restriction\gamma:\gamma\leq\beta<\omega_1\}$ are countable, 
we can find a subset $J'$ of $J$ of same cardinality and $w_\alpha$ in $\mathscr L$ extending $\tau$ such that for all $\gamma\in J',$ the set $ D'=\{w_\beta\in\mathscr L:~w_\beta~\text{extends}~w_\alpha\}\subseteq D$ has uncountable intersection with every $A_\gamma.$ However, for any $\gamma\in J'$ and for any element of $\mathscr L$ that extends $w_\alpha,$ this extension must be an element of $A_\gamma,$ since $A_\gamma$ is closed in $\mathscr L$. Thus, we have $\abs{A_\gamma\cap A_\alpha}>\aleph_0$ for $\omega_2$ many elements of $\mathscr F,$ contradicting our assumption. 
 Therefore, $\mathscr{L}$ is not a $Q$-set space.
\end{proof}

\begin{theorem}\label{Delta} If $T(o)$ is special, then $\mathscr L$ is not a $\Delta$-space.
\end{theorem}
\begin{proof} Let $T(o)$ be a countable union of antichains $A_n$. Now define $$A'_n=\{w_\beta:w_{\beta}\restriction\beta\in{A_n}\}
~\text{and}~ D_n=\mathscr L\setminus{\displaystyle\bigcup_{i<n}A'_i}.$$
Clearly, for every $n\in\omega$ $D_n\supseteq D_{n+1}$ and $\bigcap\{D_n:~{n\in\omega}\}=\emptyset$. Furthermore, consider any family $\mathscr F=\{U_n:~U_n\supseteq D_n\}$ of open sets. We show that: \[\displaystyle\bigcap_{n\in\omega}U_n\neq\emptyset\]

It suffices to show that $\mathscr L\setminus{U}_n$ is countable for each $n$. Since $\mathscr L$ is hereditarily Lindel\"{o}f, there exists a countable set $\{w_{\zeta_m}:m\in\omega\}$ such that $\bigcup_{m\in\omega}U_n(w_{\zeta_m})=U_n$, where $U_n(w_{\zeta_m})$ is a basic open neighborhood of $w_{\zeta_m}$. Fix $n$ for the moment. Let $F_m$ be the corresponding support of $U(w_{\zeta_m}),$ and let $\langle O_i^m:i<k_m\rangle$ be the finite sequence of open sets determined by the support $F_m$. We now choose an ordinal $\alpha>\max(F_m)$ $\forall~m\in\omega.$

We claim that $\forall~\gamma>\alpha,$ $w_\gamma\in{U_n},$ which will conclude that $L\setminus{U_n}$ is countable. 
To see why, we assume, for the sake of contradiction, that $w_\gamma\notin{U_n}$ for some $\gamma>\alpha$.
Hence, $\forall~m,$ $w_\gamma\restriction{F_m}\notin\prod_{i<k_m}O_i^m,$ which implies that $w_\gamma\in\bigcup_{i<n}A'_i.$ By definition, this in turn means that $w_\gamma\restriction\gamma\in\bigcup_{i<n}A_i.$ Take any extension of $w_\gamma\restriction\gamma$ in the tree $T(o)$, say $w_{\rho_0}\restriction\rho_0$. Since $\rho_0>\gamma>\alpha$ and $w_{\rho_0}\restriction\gamma=w_\gamma\restriction\gamma$, we have that $w_{\rho_0}\notin{U_n}.$ Hence, as argued above, $w_{\rho_0}\upharpoonright \rho_0\in\bigcup_{i<n}A_i.$ Now it suffices to take an infinite increasing chain of such extensions: $w_{\rho_0}\upharpoonright \rho_0\leq w_{\rho_1}\upharpoonright \rho_1\leq\cdots .$ Now by the pigeonhole principle, two of these nodes of the tree must lie in the same antichain $A_i$ for some $i<n,$ leading to a contradiction. 

Thus, $\forall\gamma>\alpha,$ we conclude that $w_\gamma\in{U_n}$ $\implies\abs{{L\setminus{U_n}}}\leq\aleph_0$. 
It now follows that $\bigcap_{n<\omega} U_n\not=\emptyset,$ which shows that $\mathscr L$ is not a $\Delta$-set space. 

\end{proof}

At this point, we raise the following question:
\begin{question} Can ${\mathscr L}$ consistently be a $\Delta$-space?
\end{question}
Clearly, by Theorem \ref{Delta}, in any model where ${\mathcal L}$ is $\Delta$, the tree $T(o)$ can't be special. In \cite{walkonordinal}, Todor\v cevi\'c proves that the Aronszajn tree defined from his $\rho_1$-function is not special in the model obtained by adding a single Cohen real. We now show the same holds true for $T(o)$. We don't know whether ${\mathscr L}$ is $\Delta$ in this model.  

\begin{proposition}
   $T(\operatorname{osc})=\{\operatorname{osc}(.,\beta)\restriction\alpha:\alpha\leq\beta<\omega_1\}$ is order isomorphic to $T(o)$, and hence, an Aronszajn tree.
   \begin{proof}
       Define the map: \begin{align*}
    \phi&:~T(o)\longrightarrow T(\operatorname{osc})\\
   &o(.,\gamma)\upharpoonright\beta\longmapsto \operatorname{osc}(.,\gamma)\upharpoonright\beta\\
\end{align*}We show that $\phi$ preserves the order $<$ induced by end extensions. However, if this is established, it would automatically show that $\phi$ is a bijective map and, hence, the trees would be order isomorphic.
To see this, fix $\alpha<\omega_1$ and consider any two countable ordinals $\gamma>\beta\geq\alpha$. \\Suppose that $o(.,\beta)\upharpoonright\alpha<o(.,\gamma)\upharpoonright\alpha$. Note that $z_\zeta$ were chosen to be rationally independent. Hence, for every $\zeta<\alpha,$ we have that: \[o(\zeta,\beta)=z_\zeta^{\operatorname{osc(\zeta,\beta)+1}}= z_\zeta^{\operatorname{osc(\zeta,\gamma)+1}}=o(\zeta,\gamma)\implies\phi(o(.,\beta)\upharpoonright\alpha)=\operatorname{osc}(\zeta,\beta)=\operatorname{osc}(\zeta,\gamma)=\phi(o(.,\gamma)\upharpoonright\alpha)\]
Therefore, $\phi(o(.,\beta)\upharpoonright\alpha)<\phi(o(.,\gamma)\upharpoonright\alpha).$
\vskip2mm
Now if $o(.,\beta)\upharpoonright \alpha\perp o(.,\gamma)\upharpoonright\alpha,$ then there would be a $\zeta<\alpha$ such that: 
$$o(\zeta,\beta)=z_\zeta^{\operatorname{osc(\zeta,\beta)+1}}\neq z_\zeta^{\operatorname{osc(\zeta,\gamma)+1}}=o(\zeta,\gamma)$$

This in turn implies that $\operatorname{osc}(\zeta,\beta)\neq\operatorname{osc}(\zeta,\beta)$ and so $\operatorname{osc}(.,\beta)\upharpoonright\alpha\perp\operatorname{osc}(.,\gamma)\upharpoonright\alpha.$
   \end{proof}
\end{proposition} 

It follows that $T(\operatorname{osc})$ is special if and only if $T(o)$ is special, so we focus on $T(\operatorname{osc})$.

\begin{lemma}\label{osc}Let $S$ be a stationary subset of $\omega_1$ and $\{\operatorname{osc}(.,\beta_\alpha)\upharpoonright\alpha:~\alpha\in{S}\}\subseteq T(\operatorname{osc})$ be given where $\beta_\alpha\geq \alpha$ for all $\alpha\in S$. Then there is a stationary subset $S'$ of $S$, such that for any pair $\alpha<\gamma$ in $S'$:
\[\operatorname{osc}(.,\alpha)\sqsubseteq\operatorname{osc}(.,\gamma)\restriction\alpha \iff ~\operatorname{osc}(.,\beta_\alpha)\restriction\alpha\sqsubseteq\operatorname{osc}(.,\beta_\gamma)\restriction\gamma\]
\end{lemma}

\begin{proof}
Without loss of generality, suppose that for every $\alpha\in S$, $\beta_\alpha>\alpha$. If it were the case that $\beta_\alpha=\alpha$ for stationarily many $\alpha,$ then the conclusion of our lemma would follow. Let $\rho_0$ be the function associated to the length of a walk defined as $\rho_0(\alpha,\beta_\alpha)=\abs{U(\alpha,\beta_\alpha)}$. By going to a stationary subset, we may assume that there is a positive integer $n$ such that for all $\alpha\in S,$ we have that $\rho_0(\alpha,\beta_\alpha)=n$. 

Now, since $C_{\beta_\alpha}\cap\alpha$ is finite for very $\alpha\in S$, we can apply a pressing down argument to obtain a stationary subset $S'$ of $S$  such that for any $\alpha\in S',$ the set $C_{\beta_\alpha}\cap\alpha$ has the same cardinality, say $s.$ 

We can now thin out our subset $S'$ of $S$ further so that for stationarily many $\alpha\in S',$ we have that $L(\alpha,\beta_\alpha)$ is the same finite set $L$. 
To establish this, first apply the pressing down lemma $s$ times to obtain a stationary set $\Gamma\subseteq S'$ such that for every $\alpha\in\Gamma$, we have that $C_{\beta_\alpha}\cap\alpha$ is the same finite set of size $s$. 
In a similar manner, as $\beta^j_\alpha$ ranging through the upper trace where $0<j<m$ (i.e., $\beta^j_\alpha\in U(\alpha,\beta_\alpha)$, and $\beta_\alpha^0=\beta_\alpha$), we can deduce that there is a stationary subset $\Gamma'$ of $\Gamma$ such that given any fixed $i<m$, we have that for every $\alpha\in\Gamma ',$ $C_{\beta^i_\alpha}\cap\alpha$ would be the same initial segment. Furthermore, we can conclude that for every $\alpha,\gamma\in\Gamma'$, $L(\alpha,\beta_\alpha)=L(\gamma,\beta_\gamma)=L$. That is because for any such $\alpha,\gamma$ and for all $j<m$, $C_{\beta^j_\alpha}\cap \alpha= C_{\beta^j_\gamma}\cap\gamma$ and the $k^{th}$ element of $L$ is defined to be: $$\operatorname{max}\big(\displaystyle\bigcup_{j=0}^k(C_{\beta^j_\alpha}\cap\alpha)\big).$$ 

Without loss of generality, we may assume that $S$ satisfies the conditions specified above. However, more refinements are still required before we can proceed to the conclusion of our lemma. First, note that for each $\alpha,$ since $\lim_{\eta\to\alpha}\min(L(\eta,\alpha))=\alpha,$ there exists a $\zeta_\alpha$ such that for all $\eta$ satisfying $\zeta_\alpha\leq\eta<\alpha,$ we have $L(\eta,\alpha)>L.$
Hence, by a pressing down argument we can find a stationary subset $S_0$ of $S$ and a countable ordinal $\zeta_0$ so that for every $\alpha\in S_0:$  
$$\forall~\eta>\zeta_0,~ L<L(\eta,\alpha)\footnote{$\max L<\min L$}$$

Since $\langle e_\alpha:\alpha<\omega_1\rangle$ (see definition \ref{max}) is a coherent family of functions, we can find a stationary subset $S_1$ of $S_0$ and an ordinal $\zeta_1\geq \zeta_0$
 such that for every $\alpha\in S_1:$ $$e_{\alpha}\upharpoonright [\zeta_1,\alpha)=e_{\beta_\alpha}\upharpoonright [\zeta_1,\alpha)$$ 

Yet again, employing the fact that $\min(L(\eta,\alpha))\to\alpha$ allows us to find a countable ordinal $\zeta_2\geq\zeta_1$ and a stationary subset $S_2$ of $S_1$ so that for every $\alpha\in S_2~\text{and}~\eta\geq\zeta_2,$ we have that: $$L(\eta,\alpha)>\zeta_1$$

Furthermore, since there are only countably many functions from $L$ to $\omega$, we can find another stationary subset $S_3$ of $S_2$ so that for every $\alpha,\gamma\in S_3,$ we have that: $$e_\alpha\restriction{L}=e_\gamma\upharpoonright L\text{ and }e_{\beta_\alpha}\restriction{L}=e_{\beta_\gamma}\upharpoonright L$$

Lastly, since the levels of $T(\operatorname{osc})$ are countable, we can refine our stationary subset $S_3$ to another stationary subset $S_4$ such that for every $\alpha<\gamma\in S_4,$ we also have that:
\begin{enumerate}
    \item $\operatorname{osc}(.,\beta_\alpha)\upharpoonright\zeta_2=\operatorname{osc}(.,\beta_\gamma)\upharpoonright\zeta_2$
    \item $\operatorname{osc}(.,\alpha)\upharpoonright\zeta_2=\operatorname{osc}(.,\gamma)\upharpoonright\zeta_2$
\end{enumerate}
 We can now present our proof. Let $n$ be the size of the set $L$ (recall $m=\rho_0(\alpha,\beta_\alpha)$ and so $n\leq m$) and consider any $\alpha<\gamma\in S_4$. 
 
We suppose $\operatorname{osc}(.,\alpha)\sqsubseteq\operatorname{osc}(.,\gamma)\restriction\alpha$ and show that $\operatorname{osc}(.,\beta_\alpha)\restriction\alpha\sqsubseteq\operatorname{osc}(.,\beta_\gamma)\restriction\gamma.$ To establish this, take any $\zeta<\alpha$. 

 If $\zeta<\zeta_2,$ then we are done since we have that $\operatorname{osc}(\zeta,\beta_\alpha)=\operatorname{osc}(\zeta,\beta_\gamma)$ by our initial refinement.  
 If $\zeta\geq\zeta_2,$ first notice that $\min(L(\zeta,\alpha))$ and $\min(L(\zeta,\gamma)$ are both bigger than $\zeta_1 > \max L.$ That in turn implies that $L(\zeta,\beta_\alpha)=L\cup{L(\zeta,\alpha)}$ and $L(\zeta,\beta_\gamma)=L\cup{L(\zeta,\gamma)}$.
 
 First, we note that since $S_4\subseteq S_3,$ we have that $e_{\beta_\alpha}\restriction{L}=e_{\beta_\gamma}\upharpoonright L,$ which implies that the number of oscillation between the pair of functions $(e_\zeta,e_{\beta_\alpha})$ on $L$ and the pair of functions $(e_\zeta,e_{\beta_\gamma})$ on $L$ are the same. 

 It remains to show that the number of oscillations between the pair $(e_\zeta,e_{\beta_\alpha})$ on $L(\zeta,\alpha)$ and the pair $(e_\zeta,e_{\beta_\gamma})$ on $L(\zeta,\gamma)$ are also the same. Note that by our hypothesis, we have that the oscillation between the pairs of functions $(e_\zeta,e_\alpha)$ and $(e_\zeta,e_\gamma)$ are the same. However, the oscillation between each pair of functions is evaluated on the lower traces $L(\zeta,\alpha)$ and $L(\zeta,\gamma),$ respectively. Now, note that the minimum of both $L(\zeta,\alpha)$ and $L(\zeta,\gamma)$ is strictly bigger than $\zeta_1,$ but we have that $e_\alpha\upharpoonright [\zeta_1,\alpha)=e_{\beta_\alpha}\upharpoonright[\zeta_1,\alpha)$. Similarly, we have that the functions $e_\gamma$ and $e_{\beta_\gamma}$ agree on $[\zeta_1,\gamma).$ Therefore, we can conclude that $\operatorname{osc}(\zeta,\beta_\alpha)\upharpoonright L(\zeta,\alpha)=\operatorname{osc}(\zeta,\beta_\gamma)\upharpoonright L(\zeta,\gamma).$

 We omit the proof of the reverse direction, as it can be proven in a similar fashion and will not be invoked in our next result.\end{proof}

Throughout the rest of the paper, given any subset $A$ of ordinals, let $\operatorname{suc}(A)$ denote the set of successor ordinals in $A$, and let $\langle f_\alpha:f_\alpha:\alpha<\omega_1\text{ a limit}\rangle$ be a fixed sequence of injective functions such that for each $\alpha$ limit, $f_\alpha:\operatorname{suc}(\alpha)\to\omega$ and for each $\alpha<\beta$ limit ordinals, $f_\alpha=^*f_\beta\upharpoonright \operatorname{suc}(\alpha)$.\footnote{for all but finitely many ordinals in $\operatorname{suc}(\alpha)$, equality is attained.} We will now obtain the non-special $T(o)$ by modifying the C-sequence using a Cohen real (see \cite{walkonordinal} for further details).

\begin{definition} Let $r\in \big([\omega]^{<\omega}\big)^\omega$ and $\langle C_\alpha:\alpha<\omega_1\rangle$ be any $C$-sequence. We define a new $C$-sequence $\langle C^r_\alpha:\alpha<\omega_1\rangle$ as follows: 
\begin{itemize}
    \item If $\alpha=\beta+1$ is a successor ordinal, then $C_\alpha^r=\{\beta\}$.
    \item If $\alpha$ is a limit ordinal, then $C_\alpha^r=C_\alpha\cap{[\zeta_\alpha^r,\alpha)}\cup\displaystyle\bigcup_{n\in\omega}{D_\alpha^r(n)},$ $$\text{where}~D_\alpha^r(n)=\{\zeta\in{\operatorname{suc}\big([C_\alpha(n),C_\alpha(n+1))\big)}\footnote{\text{let $C_\alpha(n)$ denote the $n^{th}$ element of $C_\alpha$ according to its increasing enumeration. }}:f_\alpha(\zeta)\in{r(n)}\}~\text{and}~\zeta_\alpha^r=\operatorname{sup}(\displaystyle\bigcup_{n\in\omega}{D_\alpha^r(n)})$$
\end{itemize}
  Let $\operatorname{osc}^r(.,\beta)\restriction\alpha$ denote the oscillation function using $\langle C^r_\alpha:\alpha<\omega_1\rangle$

\end{definition}

\begin{proposition}
    Let $\mathbb{P}=\operatorname{Fn}(\omega,[\omega]^{<\omega})$ the partial order adding one Cohen real. Let $\langle C_\alpha:\alpha<\omega_1\rangle$ be any $C$-sequence and $\langle f_\alpha: \alpha<\omega_1\rangle$ be any coherent sequence of functions. Then for every $\alpha<\omega_1$ and $n\in\omega$, the set \[d_{(\alpha,n)}=\{p\in\mathbb{P}:\exists~{m>n}~\text{and}~\zeta\in suc\big([C_\alpha(m),C_\alpha(m+1))\big)\ni:f_\alpha(\zeta)\in{p(m)}\}\]is a dense subset of ${\mathbb P}$. Therefore, in the Cohen extension, if $r$ is the Cohen real, then $\zeta_\alpha^r=\alpha$ and $C_\alpha^r=\displaystyle\bigcup_{n\in\omega} D_\alpha^r(n)$ for all $\alpha$.
\end{proposition}

\begin{proof}

Let $q\in\mathbb P$. We show that there is a $p\in d_{(\alpha,n)}$ such that $p\leq q.$ Note that since for any limit ordinal $\alpha$, $C_\alpha$ is a cofinal subset of $\alpha$ of order type $\omega$, we can find a positive integer $m$ where $m>\max\{\max(\operatorname{dom}(q)),n\}$ and $\min\{\operatorname{suc}\big([C_\alpha(m),C_\alpha(m+1))\big)\}\neq\emptyset.$ Let $\zeta$ be the minimum ordinal in that set. Then define $p=q^\frown\{(m,\{f_\alpha(\zeta)\})\}$, which would be the desired extension of $q$ in $d_{(\alpha,n)}.$ 

It follows from a standard density argument that $1_{\mathbb P}\Vdash\zeta_\alpha=\alpha$ and, hence, $C_\alpha^r=\displaystyle\bigcup_{n\in\omega}D_\alpha^r(n).$
\end{proof}

\begin{theorem} Let $\mathbf{M}$ be a countable transitive model of $\operatorname{ZFC}$ and $\mathbb P=\operatorname{Fn}(\omega,[\omega]^{<\omega})$. If $r$ is  a Cohen real in the forcing extension of $\mathbf{M}$, then $T(\operatorname{osc}^r)=\{\operatorname{osc}^r(.,\beta)\restriction\alpha:\alpha\leq\beta<\omega_1\}$ has no stationary antichain, and consequently, cannot be decomposed into countably many antichains. 
\end{theorem}
\begin{proof} Let $G$ be $\mathbb{P}$-generic over $\textbf{M},$ and let $r=\cup{G}$ be the Cohen real. Since any stationary set in the extension contains a ground model stationary set, let $S$ be any stationary subset of $\omega_1$ in $\textbf{M}.$  Now by Lemma \ref{osc}, it is sufficient to show that there exist ordinals $\gamma<\delta$ in $S$, such that $\operatorname{osc}^r(.,\gamma)\sqsubseteq \operatorname{osc}^r(.,\delta)\restriction\gamma$.
Given the correspondence between the dense subsets of the partial order and dense, open subsets of the model, and along with the fact that a real is Cohen generic over the ground model iff it appears in every dense open sets coded in the ground model, it suffices to show that the set: 
    $$U_S=\{x\in\big([\omega]^{<\omega}\big)^\omega:\exists~\gamma,\delta\in S\ni:\operatorname{osc}^x(.,\gamma)\sqsubseteq \operatorname{osc}^x(.,\delta)\restriction\gamma\}$$ is a dense open subset of $\big([\omega]^{<\omega}\big)^\omega$ equipped with the product topology. Note that $U_S$ is coded in the ground model.
To show that $U_S$ is dense open, it suffices to show that any $p\in\mathbb{P}$ has an extension $q$ such that the basic open set $[q]=\{g\in\big([\omega]^{<\omega}\big)^\omega:g\sqsupseteq{q}\}\subseteq{U_S}.$ To establish this, fix some $p$ and consider the minimum positive integer $n$ such that $n>\max\{\max(\operatorname{dom}(p)),\cup{\operatorname{rang}(p)}\}.$ For any $\gamma\in S,$ define the set $F_\gamma(n)=\{\zeta\in \operatorname{suc}(\gamma):f_\gamma(\zeta)\leq{n}\}$.

First, fix $S_0$ a stationary subset of $S$ and $k\in \omega$ so that for every $\alpha,\beta\in S_0,$ we have that $\abs{F_\alpha(n)}=\abs{F_\beta(n)}=k+1.$ By several applications of 
the Pressing Down Lemma, refine $S_0$ further to another stationary subset $\Gamma$ and fix $F=\{\alpha_0,\alpha_1,...,\alpha_{k}\}$ such that $\forall~\gamma\in\Gamma,$
$F_\gamma(n)=F$ and all $f_\gamma$ agree on $F$.

Without loss of generality, suppose $\alpha_{k}$ is the largest ordinal in $F$. Fix a stationary subset $\Gamma'$ of $\Gamma$ such that for every $\gamma,\gamma'\in\Gamma',$ we have that $f_\gamma\upharpoonright\alpha_{k}+1=f_{\gamma'}\upharpoonright\alpha_{k}+1.$ 

Now, fix a continuous $\epsilon$-chain of countable elementary submodels $\{N_\alpha:\alpha\in\omega_1\}$, where $N_0$ contains $\Gamma'$, $\langle\alpha_0,...,\alpha_k\rangle$, $\langle f_\alpha:\alpha\in\omega_1\rangle$, $\langle C_\alpha:\alpha\in\omega_1\rangle,$ etc.

Define $D=\{N_\alpha\cap\omega_1:\beta\leq\alpha~\text{and}~\alpha\in{\operatorname{lim}(\omega_1)}\}.$ Let $\gamma=\min(\Gamma'\cap{D})$. Consequently, $\gamma$ is of the form $N_\eta\cap\omega_1$ for some $\eta$ limit. By elementarity, our cofinal sequence $C_\gamma$ cannot be in $N_\eta$. Now, there is a $\zeta<\gamma$ such that $\abs{N_{\zeta}\cap{C_\gamma}}>n.$ Let $m+1=\abs{N_{\zeta}\cap{C_\gamma}}$.

Next, consider the set of ordinals $\alpha\in\omega_1$ for which there exists $\delta\in\Gamma' $ with the following properties: 
\begin{enumerate}
 \item${\alpha=C_\delta(m+1)}$
 \item $C_\delta\upharpoonright{m+1}=\langle C_\delta(0),\cdots,C_\delta(m)\rangle=N_{\zeta}\cap{C_\gamma}$
      \item$f_\gamma\restriction{\operatorname{suc}(C_\gamma(m)+1)}=f_\delta\restriction{\operatorname{suc}(C_\delta(m)+1})$  
         
\end{enumerate}

This set is an element of $N_{\zeta}$ since all parameters $\omega_1, \Gamma', N_{\zeta}\cap{C_\gamma}, f_\gamma\restriction{\operatorname{suc}(C_\gamma(m)+1)}$ are in $N_\zeta.$ Note that $\gamma\in\Gamma'$ and $C_\gamma(m+1)$ would serve as a witness to the fact that this set is uncountable. Therefore, there exists an $\alpha$ in this set above $\gamma+1$. 

Thus, there is a $\delta>\alpha$ such that $\alpha=C_\delta(m+1)>\gamma+1$ and the first $(m+1)$ elements of $C_\delta$ and $C_\gamma$ are the same. Additionally, $f_\gamma\restriction{\operatorname{suc}(C_\gamma(m)+1)}=f_\delta\restriction{\operatorname{suc}(C_\delta(m)+1})$ ensures that the three clauses $(1),(2),(3)$ are satisfied. 

To find an extension $q$ of $p$ in $\mathbb P$ such that the basic open set determined by $q$ is  contained in $U_S$, extend the partial function $p$ to $q$ defined on $m+1$ as follows:
$$q=\begin{cases}
     q(m)=\{f_\delta(\gamma+1)\}&\\
     q(i)=\emptyset\quad\text{if $i\notin\operatorname{dom}(p)$}
 \end{cases}$$

 Consider any extension $x\in\big([\omega]^{<\omega}\big)^{\omega}$ of $q.$ We now show that $\operatorname{osc}^x(.,\gamma)\sqsubseteq\operatorname{osc}^x(.,\delta)\upharpoonright\gamma.$ To establish this, we first need to show that $C_\delta^x\cap\gamma\sqsubseteq{C_\gamma^x}.$ The conditions given by the equalities $$f_\gamma\restriction{\operatorname{suc}(C_\gamma(m)+1)}=f_\delta\restriction{\operatorname{suc}(C_\delta(m)+1})~\text{and}~ C_\gamma\restriction{m+1}=C_\delta\restriction{m+1},$$ imply that if $$\xi\in\displaystyle\bigcup_{i=0}^{m-1}D_\delta^x(i)\subset C_\delta^x,$$
 then $\xi\in C_\gamma^x.$ Similarly, if $\xi=C_\delta(m)\in C_\delta^x.$ Hence, we only need to consider any successor ordinal $\xi$ in $(C_\delta(m),\gamma).$ If $\xi\leq\alpha_k,$ then by the choice of our finite set $F$ and the stationary set $\Gamma'$, we have that $f_\delta(\xi)=f_\gamma(\xi).$ If such $\xi$ is above $\alpha_k,$  then $f_\delta(\xi)>n>\max\{\operatorname{dom}(p),\cup\operatorname{rang}(p)\}.$ Therefore, such ordinal $\xi$ will not be included in the modified sequence $C_\delta^x$ since $x$ extends $q$ and $q$ extends $p$ and the injectivity of $f_\delta$ implies that $f_\delta(\xi)\notin x(m).$ Thus, $C_\delta^x\cap\gamma\sqsubseteq{C_\gamma^x}.$ 
 
Next, we show that $\operatorname{osc}^x(.,\gamma)\sqsubseteq \operatorname{osc}^x(.,\delta)\restriction\gamma.$ Consider any $\zeta<\gamma,$ and let us show that $\operatorname{osc}^x(\zeta,\gamma)=\operatorname{osc}^x(\zeta,\delta).$ This in turn amounts to showing that
$$e^x_\gamma=e^x_\delta\upharpoonright\gamma~\text{and}~ L^x(\zeta,\gamma)=L^x(\zeta,\delta).$$ 
To show that $e^x_\gamma=e^x_\delta\upharpoonright\gamma,$ first note that we have $\gamma+1\in{C_\delta^x}.$ That is due to the facts that $C_\delta(m)<\gamma+1<C_\delta(m+1),$ $f_\delta(\gamma+1)\in{x(m)}$ and $C_\delta^x=\cup_n{D_\delta^x(n)}.$  Therefore, for any $\xi<\gamma,$ the upper trace $U^x(\xi,\delta)$ from $\delta$ to $\xi$ with respect to the modified $C$-sequence $C_\delta^x$ would consist of 
 $\begin{cases}
       \{\gamma+1,\gamma\}\cup U^x(\xi,\gamma),\quad \text{if}~ \xi>\max (C^x_\delta\cap\gamma)\\
      U^x(\xi,\gamma),\quad \text{if otherwise.}
 \end{cases}$

Hence, invoking the definition of functions $e_\alpha,$ we can easily see that $e^x_\gamma(\xi)\leq{e^x}_\delta(\xi).$ The reverse inequality can also be established since $C_\delta^x\cap\gamma\sqsubseteq{C_\gamma^x}$ and $C_{\gamma+1}^x=\{\gamma\}$ implies that $C_{\gamma+1}^x\cap\xi=\emptyset,$ it follows that $\abs{C_\delta^x\cap\xi}\leq\abs{C_\gamma^x\cap\xi}.$ Again, by the definition of functions $e_\alpha,$ we can conclude that $e^x_\gamma(\xi)\geq{e^x}_\delta(\xi).$

Now, we show that $L^x(\zeta,\gamma)=L^x(\zeta,\delta).$
 If $\zeta\leq{\max(C_\delta^x\cap\gamma)},$ then the lower traces would be equal since $C_\delta^x\cap\gamma\sqsubseteq C_\gamma^x$ and $\min(C_\delta^x\setminus\zeta)=\min(C_\gamma^x\setminus\zeta)=\zeta$ when $\zeta={\max(C_\delta^x\cap\gamma)}.$
 If $\zeta>{\max(C_\delta^x\cap\gamma)},$ then clearly the $\min(C_\delta^x)\setminus\zeta=\gamma+1$ and since $C_{\gamma+1}^x=\{\gamma\},$ we can achieve the equality in this case as well. 
\end{proof}

\end{document}